\magnification=1200

\def\Q{{\bf {Q}}}
\def\K{{\bf {K}}}

\def\Z{{\bf Z}}

%
\catcode`@=11
%
%
\def\bibn@me{R\'ef\'erences}
\def\bibliographym@rk{\centerline{{\sc\bibn@me}}
	\sectionmark\section{\ignorespaces}{\unskip\bibn@me}
	\bigbreak\bgroup
	\ifx\ninepoint\undefined\relax\else\ninepoint\fi}
%
%
%
\let\refsp@ce=\ 
\let\bibleftm@rk=[
\let\bibrightm@rk=]
%
%
%
\def\numero{n\raise.82ex\hbox{$\fam0\scriptscriptstyle o$}~\ignorespaces}
%
%
\newcount\equationc@unt
\newcount\bibc@unt
\newif\ifref@changes\ref@changesfalse
\newif\ifpageref@changes\ref@changesfalse
\newif\ifbib@changes\bib@changesfalse
\newif\ifref@undefined\ref@undefinedfalse
\newif\ifpageref@undefined\ref@undefinedfalse
\newif\ifbib@undefined\bib@undefinedfalse
\newwrite\@auxout
%
%
\def\eqnum{\global\advance\equationc@unt by 1%
\edef\lastref{\number\equationc@unt}%
\eqno{(\lastref)}}
%
%
%
%
%
%
\def\re@dreferences#1#2{{%
	\re@dreferenceslist{#1}#2,\undefined\@@}}
\def\re@dreferenceslist#1#2,#3\@@{\def\next{#2}%
	\expandafter\ifx\csname#1@@\meaning\next\endcsname\relax
	??\immediate\write16
	{Warning, #1-reference "\next" on page \the\pageno\space
	is undefined.}%
	\global\csname#1@undefinedtrue\endcsname
	\else\csname#1@@\meaning\next\endcsname\fi
	\ifx#3\undefined\relax
	\else,\refsp@ce\re@dreferenceslist{#1}#3\@@\fi}
%
%
%
\def\newlabel#1#2{{\def\next{#1}\newl@bel#2}}
\def\newl@bel#1#2{%
	\expandafter\xdef\csname ref@@\meaning\next\endcsname{#1}%
	\expandafter\xdef\csname pageref@@\meaning\next\endcsname{#2}}
\def\label#1{{%
	\toks0={#1}\message{ref(\lastref) \the\toks0,}%
	\ignorespaces\immediate\write\@auxout%
	{\noexpand\newlabel{\the\toks0}{{\lastref}{\the\pageno}}}%
	\def\next{#1}%
	\expandafter\ifx\csname ref@@\meaning\next\endcsname\lastref%
	\else\global\ref@changestrue\fi%
	\newlabel{#1}{{\lastref}{\the\pageno}}}}
\def\ref#1{\re@dreferences{ref}{#1}}
\def\pageref#1{\re@dreferences{pageref}{#1}}
%
%
\def\bibcite#1#2{{\def\next{#1}%
	\expandafter\xdef\csname bib@@\meaning\next\endcsname{#2}}}
\def\cite#1{\bibleftm@rk\re@dreferences{bib}{#1}\bibrightm@rk}
%
%
\def\beginthebibliography#1{\bibliographym@rk
	\setbox0\hbox{\bibleftm@rk#1\bibrightm@rk\enspace}
	\parindent=\wd0
	\global\bibc@unt=0
	\def\bibitem##1{\global\advance\bibc@unt by 1
		\edef\lastref{\number\bibc@unt}
		{\toks0={##1}
		\message{bib[\lastref] \the\toks0,}%
		\immediate\write\@auxout
		{\noexpand\bibcite{\the\toks0}{\lastref}}}
		\def\next{##1}%
		\expandafter\ifx
		\csname bib@@\meaning\next\endcsname\lastref
		\else\global\bib@changestrue\fi%
		\bibcite{##1}{\lastref}
		\medbreak
		\item{\hfill\bibleftm@rk\lastref\bibrightm@rk}%
		}
	}
\def\endthebibliography{\egroup\par}
%
%
%
\def\@closeaux{\closeout\@auxout
	\ifref@changes\immediate\write16
	{Warning, changes in references.}\fi
	\ifpageref@changes\immediate\write16
	{Warning, changes in page references.}\fi
	\ifbib@changes\immediate\write16
	{Warning, changes in bibliography.}\fi
	\ifref@undefined\immediate\write16
	{Warning, references undefined.}\fi
	\ifpageref@undefined\immediate\write16
	{Warning, page references undefined.}\fi
	\ifbib@undefined\immediate\write16
	{Warning, citations undefined.}\fi}
%
%
\immediate\openin\@auxout=\jobname.aux
\ifeof\@auxout \immediate\write16
  {Creating file \jobname.aux}
\immediate\closein\@auxout
\immediate\openout\@auxout=\jobname.aux
\immediate\write\@auxout {\relax}%
\immediate\closeout\@auxout
\else\immediate\closein\@auxout\fi
%
%
\input\jobname.aux
\immediate\openout\@auxout=\jobname.aux
%
%
\catcode`@=12

%
\catcode`@=11
\def\bibliographym@rk{\bgroup}
%
%
\outer\def\bye{ 	\par\vfill\supereject\end}

   \def\resp{{\it resp.,  }}

\def\house#1{\setbox1=\hbox{$\,#1\,$}%
\dimen1=\ht1 \advance\dimen1 by 2pt \dimen2=\dp1 \advance\dimen2 by 2pt
\setbox1=\hbox{\vrule height\dimen1 depth\dimen2\box1\vrule}%
\setbox1=\vbox{\hrule\box1}%
\advance\dimen1 by .4pt \ht1=\dimen1
\advance\dimen2 by .4pt \dp1=\dimen2 \box1\relax}

  \def\eps{{\varepsilon}}

 \def\ens{\enspace} \def\noi{\noindent}

\def\build#1_#2^#3{\mathrel{\mathop{\kern 0pt#1}\limits_{#2}^{#3}}}

\def\date {le\ {\the\day}\ \ifcase\month\or janvier
\or fevrier\or mars\or avril\or mai\or juin\or juillet\or
ao\^ut\or septembre\or octobre\or novembre
\or d\'ecembre\fi\ {\oldstyle\the\year}}

\font\fivegoth=eufm5 \font\sevengoth=eufm7 \font\tengoth=eufm10

\newfam\gothfam \scriptscriptfont\gothfam=\fivegoth
\textfont\gothfam=\tengoth \scriptfont\gothfam=\sevengoth

\def\smallsquare{\vbox{\hrule\hbox{\vrule height 1 ex\kern 1 ex\vrule}\hrule}}

\def\og{\leavevmode\raise.3ex\hbox{$\scriptscriptstyle 
\langle\!\langle\,$}}
\def \fg {\leavevmode\raise.3ex\hbox{$\scriptscriptstyle 
\!\rangle\!\rangle\,\,$}}

\def\rme{{\rm e}}

\def\mueff{{\mu_{\rm eff}}}

\centerline{}

\vskip 8mm

\centerline{\bf Effective simultaneous rational approximation} 

\smallskip

\centerline{\bf to pairs of real quadratic numbers}

\vskip 9mm

\centerline{Y{\sevenrm ANN} B{\sevenrm UGEAUD} \footnote{}{\rm
2010 {\it Mathematics Subject Classification : } 11J13; 11D09, 11J86.}}

\vskip 8mm

\hfill{\it To the memory of Naum Ilich Feldman (1918--1994)}

{\narrower\narrower
\vskip 12mm

\proclaim Abstract. {
Let $\xi, \zeta$ be quadratic real numbers in distinct quadratic fields. 
We establish the existence of effectively computable, 
positive real numbers $\tau$ and $c$, such that, for every integer $q$ with $q > c$ we have
$$
\max\{\|q \xi \|, \|q \zeta\| \} > q^{-1 + \tau},
$$ 
where $\| \cdot \|$ denotes the distance to the nearest integer. 
}

}

\vskip 10mm

\centerline{\bf 1. Introduction and results}

\vskip 5mm

Let $\xi$ be an irrational real number. The real number $\mu$ 
is an irrationality measure for $\xi$ if there exists 
a positive real number $c(\xi)$ such that every rational number ${p \over q}$
with $q \ge 1$ satisfies 
$$
\Bigl| \xi - {p \over q} \Bigr| > {c(\xi) \over q^{\mu}}.
$$
If, moreover, the constant $c(\xi)$ is effectively computable, 
then $\mu$ is an 
effective irrationality measure for $\xi$. 
We denote by $\mu (\xi)$ (\resp  $\mueff (\xi)$) the infimum of the 
irrationality measures (\resp  effective irrationality measures) for $\xi$ and call it 
the irrationality exponent  
(\resp effective irrationality exponent) of $\xi$.  
It follows from the theory of continued fractions that $\mu(\xi) \ge 2$ and an easy 
covering argument shows that equality holds for almost all $\xi$, with respect to the 
Lebesgue measure. Furthermore, if $\xi$ is real algebraic of degree $d \ge 2$, 
then Liouville's inequality implies that $\mueff (\xi) \le d$, while 
Roth's theorem asserts that $\mu (\xi) = 2$. 
To get better upper bounds for 
the effective irrationality exponents of algebraic numbers is a 
notorious challenging problem. 

The first result of this type 
was obtained in 1964 by Alan Baker \cite{Ba64}, who 
established that $\mueff (\root 3\of 2) \le 2.955$, but his method applies 
only to a very restricted class of algebraic numbers. 
A few years later, in 1971, Feldman \cite{Fe71}, by means of a refinement of the lower bounds 
for linear forms in logarithms of algebraic numbers established by Baker, proved that the 
effective irrationality exponent  
of an arbitrary real algebraic number of 
degree greater than two is strictly less than its degree; see also \cite{BiBu00} for a 
proof depending on lower bounds 
for linear forms in only two logarithms. 
Subsequently, Bombieri \cite{Bo93,Bo02}  
gave in 1993 an alternative proof of Feldman's result, 
completely independent of the theory of linear forms in logarithms and based on the
Thue--Siegel Principle. 
Further results and bibliographic references can be found in \cite{BuLiv18}, see 
in particular Section 4.10. 

In this note, we are concerned with the simultaneous approximation to pairs of real numbers
by rational numbers having the same denominator. We extend the above 
definition of (effective) irrationality exponent as follows. 
Let $\xi, \zeta$ be real numbers such that $1, \xi, \zeta$ are linearly independent 
over the rational numbers. The real number $\mu$ 
is a simultaneous irrationality measure for the pair $(\xi, \zeta)$ if there exists 
a positive real number $c(\xi, \zeta)$ such that, for every integer triple $(p, q, r)$
with $q \ge 1$, we have 
$$
\max \Bigl\{ \Bigl| \xi - {p \over q} \Bigr|,  \Bigl| \zeta - {r \over q} \Bigr|  \Bigr\} > 
{c(\xi, \zeta) \over q^{\mu}}.
$$
If, moreover, the constant $c(\xi, \zeta)$ is effectively computable, 
then $\mu$ is an effective irrationality measure for the pair $(\xi, \zeta)$.
We denote by $\mu (\xi, \zeta)$ (\resp  $\mueff (\xi, \zeta)$) the infimum of the 
irrationality measures (\resp  effective irrationality measures) for the pair $(\xi, \zeta)$ and call it 
the irrationality exponent  
(\resp effective irrationality exponent) of the pair $(\xi, \zeta)$.

Let $\xi, \zeta$ be real numbers such that $1, \xi, \zeta$ are linearly independent 
over the rational numbers. An easy application of Minkowski's theorem implies 
that $\mu (\xi, \zeta) \ge {3 \over 2}$ and a covering lemma shows that equality holds 
for almost all pairs $(\xi, \zeta)$, with respect to the planar Lebesgue measure. 
Schmidt \cite{Schm67} established 
that $\mu (\xi, \zeta) = {3 \over 2}$ if $\xi$ and $\zeta$ are both real and algebraic. 
His result is ineffective and gives no better information on 
$\mueff (\xi, \zeta)$ than the obvious inequality 
$$
\mueff (\xi, \zeta) \le \max\{ \mueff(\xi), \mueff(\zeta) \}.
$$
The particular case where $\xi$ and $\zeta$ are quadratic numbers in distinct 
number fields is of special interest. The obvious upper bound 
$\mueff (\xi, \zeta) \le 2$ has been improved in some cases, 
in particular by Rickert \cite{Ri93} (see his paper for earlier references), who 
established among other results that 
$$
\max \Bigl\{ \Bigl| \sqrt{2} - {p \over q} \Bigr|,  \Bigl| \sqrt{3} - {r \over q} \Bigr| \Bigr\} > 
{10^{-7} \over q^{1.913}}, \quad \hbox{for integers $p, q, r \ge 1$},    
$$
and subsequently by Bennett \cite{Be95,Be96}. 
The method used in these papers applies only to a very restricted 
class of pairs $(\xi, \zeta)$ of quadratic numbers. 

The purpose of the present note is to show how the theory of linear forms 
in logarithms (or, alternatively, Bombieri's method) allows us to improve the trivial upper bound 
$\mueff (\xi, \zeta) \le 2$ for all quadratic real numbers $\xi$ and $\zeta$ 
in distinct quadratic fields.

\proclaim Theorem 1.1. 
Let $\xi, \zeta$ be real quadratic numbers in distinct quadratic fields. 
Let $R_\xi$ and $R_\zeta$ denote the regulators
of the fields $\Q(\xi)$ and $\Q(\zeta)$, respectively. 
Then, there exists an absolute, positive, effectively computable real number $c_1$ such that 
$$
\mueff (\xi, \zeta) \le 2 - (c_1 R_\xi R_\zeta )^{-1}.    \eqno (1.1)
$$
In particular, if $a, b$ are positive integers such that none of $a, b$, 
and $ab$ is a perfect square, then 
there exists an absolute, positive, effectively computable real number $c_2$ such that 
$$
\mueff (\sqrt a, \sqrt b) \le 2 - (c_2 \sqrt{ab} (\log a) (\log b) )^{-1}. 
$$

The last assertion of Theorem 1.1 is an immediate consequence of the 
first one, since for any square-free integer $D \ge 2$ the regulator $R_D$ of the 
quadratic field generated by $\sqrt{D}$ satisfies
$$
R_D < \sqrt{D} (1 + \log \sqrt{D}),   \eqno (1.2) 
$$
see e.g. \cite{JaLuWi95}.

Theorem 1.1 is by no means surprising. It is ultimately a consequence of the 
quantity $B'$, which has its origin in Feldman's papers \cite{Fe68,Fe71} and is the key tool 
for his effective improvement of Liouville's bound; see Theorem 2.1 and the discussion below it.  
Other consequences of the quantity $B'$ can be found in \cite{BuLiv18} and 
in the recent papers \cite{Bu18,BuEv17,BuEvGy18}. 

We present a proof of Theorem 1.1 together with a proof of a slightly weaker version of it, with 
$R_\xi R_\zeta$ replaced by $R_\xi R_\zeta \log (R_\xi R_\zeta)$ in (1.1). 
For the latter result, we apply an estimate for 
linear forms in three logarithms, while the former is derived from a result of 
Bombieri \cite{Bo93} (and can also be derived from an estimate for 
linear forms in only two logarithms). This is in accordance with the improvements on 
Liouville's bound obtained by these two methods. Namely, for an algebraic number $\xi$ 
of degree $d$ at least equal to $3$, denoting by $R_\xi$ the regulator of the number field 
generated by $\xi$, it follows from the theory of linear forms in logarithms 
and from Bombieri's method, respectively, 
that there exist effectively computable, positive real numbers $c_3$ and $c_4$ such that 
$$
\mueff (\xi) \le d - (c_3 R_\xi   \log R_\xi  )^{-1} 
$$
and
$$
\mueff (\xi) \le d - (c_4 R_\xi)^{-1},
$$
respectively; see e.g. \cite{Bu98}.

The last assertion of Theorem 1.1 is equivalent to the following statement on 
systems of Pellian equations. 

\proclaim Theorem 1.2. 
Let $a, b$ be positive integers such that none of $a, b$, and $ab$ is a perfect square. 
Let $u, v$ be non-zero integers. 
Then, there exists an effectively computable, absolute real number $c_5$ such that 
all the solutions in positive integers $x, y, z$ of the system of Pellian equations
$$
x^2 - a y^2 = u, \quad z^2 - b y^2 = v 
$$
satisfy
$$
\max\{ x, y, z \} \le (\max\{ |u|, |v|, 2 \})^{c_5 \sqrt{ab} (\log a) (\log b) }. 
$$




\vskip 5mm

\centerline{\bf 2. Auxiliary results}

\vskip 5mm

As usual, $h(\alpha)$ denotes the (logarithmic) Weil height of the 
algebraic number $\alpha$. 
Our auxiliary result for the proof of (a slightly weaker version of) Theorems 1.1 and 1.2 
is a particular case of Theorem 2.1 of \cite{BuLiv18}, which 
essentially reproduces a theorem of Waldschmidt \cite{Wa93,WaLiv}. 

\proclaim Theorem 2.1. 
Let $n \ge 1$ be an integer. 
Let $\alpha_1, \ldots, \alpha_n$ be non-zero algebraic numbers. 
Let $b_1, \ldots , b_n$ be integers with $b_n \not= 0$.
Let $D$ be the degree over $\Q$ 
of the number field $\Q(\alpha_1, \ldots, \alpha_n)$. 
Let $A_1, \ldots, A_n$ be real numbers with
$$
\log A_j \ge \max
\Bigl\{ h(\alpha_j), { \rme \over D} |\log \alpha_j|, {1 \over D} \Bigr\},
\qquad 1\le j \le n.
$$
Let $B'$ be a real number satisfying 
$$
B'\ge 3D, \quad 
B' \ge \max_{1 \le j \le n-1} 
\Bigl\{ {|b_n| \over \log A_j} + {|b_j| \over \log A_n} \Bigr\}. 
$$
If $b_1 \log \alpha_1 + \cdots + b_n \log \alpha_n$ is nonzero, then we have
$$
\eqalign{
\log |b_1 \log \alpha_1 & + \cdots + b_n \log \alpha_n|   \cr
& \ge - 2^{n + 26} \, n^{3n + 9} \, D^{n+2} \, \log (3D) \, \log A_1 \ldots \log A_n \log B'. \cr} 
$$

The quantity $B'$ in Theorem 2.1, which replaces the quantity 
$$
B = \max\{3D, |b_1|, \ldots , |b_n|\}
$$ 
occurring in earlier estimates of Baker, originates in Feldman's papers \cite{Fe68,Fe71}. It is 
a consequence of the use of the functions $x \mapsto {x \choose k}$ instead of 
$x \mapsto x^k$ in the construction of the auxiliary function. 
The key point is the presence of the factor $\log A_n$ in the denominator in the
definition of $B'$. It is of great interest when $b_n = 1$ and $\log A_n$ is
large, since it then allows us, roughly 
speaking, to replace $B$ by $B / (\log A_n)$.

The auxiliary result for the proof of Theorems 1.1 and 1.2 
is a particular case of Theorem 2 of Bombieri \cite{Bo93}. 
Actually, since the dependence in the parameters
$d$ and $\kappa$ occurring in this theorem has been improved in \cite{Bu98}, 
we choose to quote below a particular case of Th\'eor\`eme 1 of \cite{Bu98}.

\proclaim Theorem 2.2.
Let $\K$ be a real number field of degree $d$.
Let $\Gamma$ be a finitely generated subgroup of  $\K^{*}$ and consider a 
system $\xi_1, \ldots, \xi_t$ of generators of $\Gamma / {\rm tors}$.
Let $\xi$ in $\Gamma$, $A$ in $\K^{*}$ and $\kappa > 0$ be such that $\kappa \le 1$ and
$$
0 < \vert 1 - A \xi \vert <  {\rm e}^{-\kappa h(A \xi)} <1.
$$
Setting
$$
C = 4.10^{19} \, d^4 \, {(\log 3d)^7 \over \kappa} \, \log^*{d \over \kappa},
\qquad Q = ( 2 t C)^t \, \prod_{i=1}^t  h(\xi_i),
$$
we have the upper bound 
$$
h(\xi) \le 10 \,  Q \max \{h(A), Q \}.
$$

Bombieri's original proof of Theorem 2.2 (upto the dependence on $d$ and $\kappa$) 
is independent 
of the theory of linear forms in logarithms. An alternative proof, given in \cite{Bu98}, 
depends on lower estimates for linear forms in two logarithms (a careful reader can 
observe that, while the proof of Th\'eor\`eme 1 of \cite{Bu98} rests on estimates for linear forms 
in three logarithms, estimates for linear forms in two logarithms are enough 
to establish Theorem 2.2 above, and even with a better numerical constant, since 
we have assumed that $\K$ is a real number field) combined with a lemma of 
geometry of numbers from \cite{Bo93}. To deduce Theorem 2.2 from 
estimates for linear forms in two logarithms, the crucial 
ingredient is ultimately the presence of the factor $B'$ in these estimates.

\vskip 5mm

\centerline{\bf 3. Proofs}

\vskip 5mm

We start with the proof of (a slightly weaker version of) Theorem 1.2. 
Let $a, b$ be positive integers such that $1, \sqrt a, \sqrt b$ 
are linearly independent over the rationals. 
Let $u, v$ be nonzero integers and consider the system of Pellian equations
$$
x^2 - a y^2 = u, \quad z^2 - b y^2 = v, \quad \hbox{in positive integers $x, y, z$}.  \eqno (3.1) 
$$
Set 
$$
U = \max\{|u|, |v|, 2 \} \quad \hbox{and} \quad X = \max\{x, y, z\}. 
$$
It is well-known \cite{BaDa69,Pinch88} that the theory of linear forms in logarithms allows us to 
bound effectively $X$ in terms of $U$. Our goal is to show that we can get a bound which is 
polynomial in $U$.

Let $\eps$ and
$\eta$ be the fundamental totally positive units of the 
rings of integers of the fields $\Q(\sqrt a)$ and $\Q(\sqrt b)$, respectively, normalized 
to be greater than $1$. We note that $\xi$ and $\eta$ are at least 
equal to $(1 + \sqrt{5})/2$. 

Let $x$, $y$, and $z$ be positive integers satisfying (3.1). 
Since the norm over $\Q$ of $x + y\sqrt a$ (\resp $z + y \sqrt b$) is $u$ (\resp $v$),
there exist nonnegative integers $m, n$ and algebraic numbers 
$\alpha$ in $\Q(\sqrt a)$ and $\beta$ in $\Q( \sqrt b)$
such that 
$$
\alpha \ge |\alpha^\sigma|, \ens \beta \ge |\beta^\sigma|, \ens
\alpha \eps^{-1} \le |\alpha^\sigma| \eps, \ens 
\beta \eta^{-1} \le |\beta^\sigma| \eta,   \eqno (3.2) 
$$
$$
x + y\sqrt a = \alpha \eps^m, \quad {\rm and} \quad
z + y \sqrt b = \beta \eta^n, 
$$
where the superscript ${ \cdot }^{\sigma}$ denotes the Galois conjugacy. 

Since $\eps^\sigma= \eps^{-1}$ and $\eta^\sigma = \eta^{-1}$, we have 
$$
2 y\sqrt a = \alpha \eps^m - \alpha^{\sigma} \eps^{-m}
$$
and
$$
2 y\sqrt b =  \beta \eta^n - \beta^{\sigma} \eta^{-n}. 
$$
Set
$$
\Lambda = |\alpha \beta^{-1} \sqrt{b/a} \, \eps^m \eta^{-n} - 1| 
= |\alpha^\sigma \beta^{-1} \sqrt{b/a} \, \eps^{-m} \eta^{-n} - \beta^\sigma \beta^{-1} \eta^{-2n}|.
\eqno (3.3)
$$
Clearly, $\Lambda$ is nonzero.

Set
$$
U_0 = \max\{U, ab, \eps^2, \eta^2\}   \eqno (3.4)
$$
Observe that $\alpha = |u| / |\alpha^\sigma|$, $\beta = |v| / |\beta^\sigma|$, (3.2), 
and (3.4) imply that 
$$
\alpha^2 \le |u| \eps^2 \le U_0^2, \quad 
\beta^2 \le |v| \eta^2 \le U_0^2,      \eqno (3.5)
$$
and
$$
\eqalign{
h(\alpha \beta^{-1} \sqrt{b/a}) & \le h(\alpha) + h(\beta) + h(\sqrt a) + h(\sqrt b) \cr
& \le \log \alpha + \log \beta + (\log a) / 2 + (\log b) / 2  \le 3 \log U_0. \cr}
$$

Assume first that 
$$
\max\{m \log \eps, n \log \eta\} \ge 12 \log U_0.      \eqno (3.6) 
$$
Observe that (3.3), (3.4), and (3.5) imply that
$$
\log \Lambda \le -  n \log \eta + 2 \log U_0,      \eqno (3.7)
$$
and
$$
|m \log \eps - n \log \eta| \le 4  \log U_0,   \eqno (3.8)
$$
thus, by (3.6), we get 
$$
\log \Lambda \le - \max\{m \log \eps, n \log \eta\} + 6 \log U_0 
\le - {\max\{m \log \eps, n \log \eta\} \over 2}.   \eqno (3.9)
$$
It then follows from Theorem 2.1 applied with $\alpha_1 = \eps$, 
$\alpha_2 = \eta$, $\alpha_3 = \alpha \beta^{-1} \sqrt{b/a}$ that 
$$
\log \Lambda \gg - (\log U_0) \, (\log \eps) \, (\log \eta) \, 
\log^* {\max\{m, n\} \over  \log U_0  },     \eqno (3.10) 
$$
where we write $\log^*$ for the function $\max\{1, \log \}$. 
Here and below, the numerical constant implied by $\ll$ is positive, absolute, and 
effectively computable. 

The combination of (3.9) with (3.10) gives
$$
\max\{m \log \eps, n \log \eta\}  \ll   (\log U_0) \, (\log \eps) \, (\log \eta) \, 
\log^* {\max\{m, n\} \over  \log U_0  }. 
$$
We deduce that 
$$
X \ll \max\{m \log \eps, n \log \eta\} \ll  (\log \eps) \, (\log \eta) \, 
\log^* (\max\{ \log \eps, \log \eta\} ) \, \log U_0, 
$$
while $X \ll \log U_0$ if (3.6) is not satisfied.  

Consequently, no matter if (3.6) holds or not, 
there exist an effectively computable positive real number $C_1$, depending only 
on $a$ and $b$, and an effectively computable positive, absolute real number $c_6$ such that 
$$
X \le C_1 \, U^{c_6 ( \log \eps ) \, ( \log \eta )  \, 
\log^* (\max\{ \log \eps, \log \eta\} ) }.     \eqno (3.11) 
$$
Combined with the upper bound (1.2), this gives Theorem 1.2 upto an extra 
logarithmic factor. 

For the proof of (a slightly weaker version of) Theorem 1.1, 
without any loss of generality, we may assume that $\xi, \eta$ are positive integers $a, b$ 
as above. Then, keeping our notation, it follows from (3.11) that there exists 
an effectively computable positive real number $C_2$, depending only 
on $a$ and $b$, such that 
$$
\eqalign{
\max \Bigl\{ \Bigl| \sqrt a  - {x \over y} \Bigr|,  \Bigl| \sqrt b - {z \over y} \Bigr|  \Bigr\}  
& = {1 \over y^2} \, \max\{|x^2 - a y^2|, |z^2 - b y^2| \} \cr
& \ge {1 \over y^2} \, \Bigl( {X \over C_1} \Bigr)^{1 / (c_6 ( \log \eps ) \, ( \log \eta )  \, 
\log^* (\max\{ \log \eps, \log \eta\} )) } \cr 
& \ge {C_2 \over y^{2 - 1 / (c_6 ( \log \eps ) \, ( \log \eta )  \, 
\log^* (\max\{ \log \eps, \log \eta\} )) }}. 
\cr}
$$
Combined with (1.2), this completes the proof of Theorem 1.1 upto an extra 
logarithmic factor.

\vskip 3mm

It remains for us to explain how to deduce Theorems 1.1 and 1.2 from 
Theorem 2.2, applied with $\Gamma$ being the subgroup generated by $\eps$ and $\eta$,
$$
A = \alpha \beta^{-1} \sqrt{b/a}, \enspace \xi_1 = \eps, \enspace 
\xi_2 = \eta, \quad \hbox{and} \quad 
\xi = \eps^m \eta^{-n}. 
$$
Note that 
$$
h(A \xi) \le h(A) + m \log \eps + n \log \eta 
\le 3 \log U_0  + m \log \eps + n \log \eta.   \eqno (3.12)
$$
Assume that (3.6) holds. 
By combining (3.6), (3.7), (3.8), and (3.12) we get 
$$
\log \Lambda  \ll  
- \log U_0  - m \log \eps - n \log \eta  \ll - h(A \xi).
$$
It then follows from Theorem 2.2 that  
$$
h(\xi) \ll  \bigl( (\log \eps) \, (\log \eta)  h(A) + (\log \eps)^2 \, (\log \eta)^2 \bigr).
$$
Since $h(A) \le 3 \log U_0$ and 
$$
X \ll \max\{m \log \eps, n \log \eta\} \le 4 h(\xi), 
$$
there exist an effectively computable positive real number $C_3$, depending only 
on $a$ and $b$, and an effectively computable positive, absolute real number $c_7$ such that 
$$
X \le C_3 \, U^{c_7 ( \log \eps ) \, ( \log \eta )  }.     \eqno (3.13)
$$
By increasing $c_7$ and $C_3$ if necessary, we see that (3.13) also holds if (3.6) is not satisfied. 
Then, proceeding as below (3.11), we establish Theorems 1.1 and 1.2.

\vskip 6mm

\noindent 
{\bf Acknowledgements. }
The idea of this note came immediately at the end of the workshop 
co-organized by Andrej Dujella in Dubrovnik at the 
end of June 2019. I am very pleased to thank him and the speakers of 
Friday morning.

\vskip 18mm

\centerline{\bf References}

\vskip 7mm

\beginthebibliography{999}

\medskip

\bibitem{Ba64}
A. Baker,
{\it Rational approximations to $\root 3\of 2$ and other algebraic numbers}, 
Quart. J. Math. Oxford Ser. 15 (1964), 375--383. 

\bibitem{BaDa69}
A. Baker and H. Davenport, 
{\it The equations $3 x^2 - 2 = y^2$ and $8 x^2 - 7= z^2$}, 
Quart. J. Math. Oxford Ser. (2) 20 (1969), 129--137.

\bibitem{Be95}
M. A. Bennett,  
{\it Simultaneous approximation to pairs of algebraic numbers}. 
In: Number theory (Halifax, NS, 1994), 55--65,
CMS Conf. Proc., 15, Amer. Math. Soc., Providence, RI, 1995.

\bibitem{Be96}
M. A. Bennett,  
{\it Simultaneous rational approximation to binomial functions}, 
Trans. Amer. Math. Soc. 348 (1996), 1717--1738.

\bibitem{BiBu00}
Yu. Bilu et Y. Bugeaud,
{\it D\'emonstration du th\'eor\`eme de
Baker-Feldman via les formes lin\'eaires en deux logarithmes}, 
J. Th\'eor. Nombres Bordeaux 12 (2000), 13--23. 

\bibitem{Bo93}
E. Bombieri,
{\it Effective Diophantine approximation on ${\bf G}_m$}, 
Ann. Scuola Norm. Sup. Pisa Cl. Sci. 20 (1993), 61--89.

\bibitem{Bo02}
E. Bombieri,
{\it Forty years of effective results in Diophantine theory}. 
In: A panorama of number theory or the view from Baker's garden (Z\"urich, 1999), 194--213, 
Cambridge Univ. Press, Cambridge, 2002.

\bibitem{Bu98}
Y. Bugeaud, 
{\it Bornes effectives pour les solutions des \'equations en
$S$-unit\'es et des \'equations de Thue-Mahler},
J. Number Theory 71 (1998), 227--244.

\bibitem{Bu18}
Y. Bugeaud, 
{\it On the digital representation of integers with bounded prime factors}, 
Osaka J. Math. 55 (2018), 315--324. 

\bibitem{BuLiv18}
Y. Bugeaud, 
Linear forms in logarithms and applications. 
IRMA Lectures in Mathematics and Theoretical Physics, 28. 
European Mathematical Society (EMS), Z\"urich, 2018.

\bibitem{BuEv17}
Y. Bugeaud and J.-H. Evertse, 
{\it $S$-parts of terms of integer linear recurrence sequences}, 
Mathematika 63 (2017), 840--851.

\bibitem{BuEvGy18}
Y. Bugeaud, J.-H. Evertse, and K. Gy\H ory, 
{\it $S$-parts of values of of univariate polynomials, binary forms 
and decomposable forms at integral points}, 
Acta Arith. 184 (2018), 151--185.

\bibitem{Fe68}
{N. I. ~Fel'dman},
{\it Improved estimate for a linear form of the logarithms of 
algebraic numbers}, 
{Mat. Sb.} {  77} (1968), 393--406 (in Russian);
English translation in Math. USSR. Sb.  6 (1968) 423--436.

\bibitem{Fe71}
{N. I. ~Fel'dman},
{\it An effective refinement of the exponent in Liouville's theorem}, 
Iz. Akad. Nauk SSSR, Ser. Mat. 35 (1971), 973--990 (in Russian);
English translation in Math. USSR. Izv. 5 (1971) 985--1002.

\bibitem{JaLuWi95}
M. J. Jacobson, Jr., R. F. Lukes, and H. C. Williams, 
{\it An investigation of bounds for the regulator of quadratic fields}, 
Experiment. Math. 4 (1995), 211--225.

\bibitem{Pinch88}
R. G. E. Pinch, 
{\it Simultaneous Pellian equations}, 
Math. Proc. Cambridge Philos. Soc. 103 (1988), 35--46.

\bibitem{Ri93}
J. H. Rickert, 
{\it Simultaneous rational approximations and related Diophantine equations}, 
Math. Proc. Cambridge Philos. Soc. 113 (1993),  461--472.

\bibitem{Schm67}
{W. M. Schmidt},
{\it On simultaneous approximations of two algebraic numbers by rationals}, 
Acta Math. 119 (1967), 27--50.

\bibitem{Wa93}
M. Waldschmidt,
{\it Minorations de combinaisons 
lin\'eaires de logarithmes de nombres
alg\'ebri\-ques}, Canadian J. Math. 45 (1993), 176--224.

\bibitem{WaLiv}
M. Waldschmidt,
Diophantine Approximation on Linear Algebraic Groups. 
Transcendence Properties of the Exponential Function in Several Variables, 
Grundlehren Math. Wiss. 326, Springer, Berlin, 2000.

\vskip 1cm

\noi Yann Bugeaud

\noi Institut de Recherche Math\'ematique Avanc\'ee, U.M.R. 7501

\noi Universit\'e de Strasbourg et C.N.R.S.

\noi 7, rue Ren\'e Descartes

\noi 67084 STRASBOURG \ \ (France)

\vskip 2mm

\noi e-mail : {\tt bugeaud@math.unistra.fr}

\bye